\documentclass{article}
\usepackage[utf8]{inputenc}
\usepackage{xcolor}
\usepackage[normalem]{ulem}
\usepackage{tikz}
\usepackage[english]{babel}
\usepackage[T1]{fontenc}
\usepackage{physics}
\usepackage{mathtools}
\usepackage{amsthm}
\usepackage{amssymb}
\usepackage{fancybox}
\usepackage[letterpaper,top=2cm,bottom=2cm,left=3cm,right=3cm,marginparwidth=1.75cm]{geometry}

\usepackage{amsmath}
\usepackage{graphicx}
\usepackage[colorlinks=true, allcolors=blue]{hyperref}

\theoremstyle{plain}
\newtheorem{theorem}{Theorem}[section]
\newtheorem{lemma}[theorem]{Lemma}
\newtheorem{corollary}[theorem]{Corollary}

\newtheorem{conj}[theorem]{Conjecture}

\theoremstyle{definition}

\theoremstyle{remark}

\newenvironment{Proof}{
	\noindent{
		\bf {Proof:}
		}
	}{$\hfill \blacksquare $\newline}

\begin{document}

\title{Group irregularity strength of disconnected graphs}

\author{Sylwia Cichacz, Barbara Krupińska\footnote{This work was partially supported by the Faculty of Applied Mathematics AGH UST statutory tasks within subsidy of Ministry of Science and Higher Education.}
\normalsize \\AGH University of Krakow, \vspace{2mm} Poland\\
}
\maketitle
\begin{abstract}
We investigate the \textit{group irregular strength} $(s_g(G))$ of graphs, i.e the smallest value of $s$ such that for any Abelian group $\Gamma$ of order $s$ exists a function $g\colon E(G) \rightarrow \Gamma$ such that sums of edge labels at every vertex is distinct. We give results for bound and exact values of $(s_g(G))$ for graphs without small stars as components. 
\end{abstract}

\section{Introduction}
One of the first facts we can observe in any simple graph $G$ is that there are at least two vertices with equal amounts of adjacent edges (of the same degree).
However, if we consider edge labeling $f\colon E(G) \rightarrow \{1\ldots s\}$ and calculate the weight of each vertex $v$ as the sum of labels on the adjacent edges, then it is possible that every vertex could have a different weight.
A labeling $f$ that produces a unique weighted degree for every vertex $v\in V(G)$ is called \textit{irregular}.
The maximum label $s$ is called the \textit{strength} of the labeling.
\textit{Irregularity strength of $G$}, denoted by $s(G)$, is the minimum value of $s$ allowing some irregular labeling.\\
Chartrand \cite{IrrNet} was the first to introduce the problem of finding $s(G)$ and many authors investigated this problem further.
In general it is known that $s(G) \leq n -1$ for any graph
$G$ of order $n$  containing no isolated edges and at most one isolated vertex, with the exception of the graph $K_3$, as proved by Aigner
and Triesch \cite{AigTrie} and Nierhoff \cite{TigBouIrrStr}. Though the family of stars witnesses the tightness of this upper
bound, this result can be improved for graphs with large enough minimum degree $\delta$.
The best general upper bound was proven by Kalkowski \cite{UppBouIrrStr} and states that $s(G)\leq 6n/\delta +6$.
This result  was improved recently by  Przybyło and Wei, who proved that for any $\varepsilon \in  (0, 0.25)$
there exist absolute constants $c_1, c_2$ such that for all graphs $G$ on $n$ vertices with minimum degree
$\delta\geq 1$ and without isolated edges, $s(G) \leq  n/\delta(1 + c_1/\delta^{\varepsilon}) + c_2$ \cite{PrzWei}.


Assume $\Gamma$ is an Abelian group of order $m \geq n$ with the operation denoted by $+$ and the identity element by $0$.
To denote $a + a + \ldots + a$, where element $a$ appears $k$ times, we will use $ka$, we will write $-a$ to denote the inverse of $a$, and instead of $a+(-b)$ we will use $a-b$.
Furthermore the notation $\sum\limits_{a\in S}a$ will be used as a shortened version of $a_1+a_2+a_3+\ldots$, where $a_1,a_2,a_3,\ldots$ are elements of the $S$ set.\\

In this paper we will focus on edge labeling $f:E(G)\rightarrow \Gamma$ which implies the weighted degrees defined as the sum in $\Gamma$:
$$w(v) = \sum\limits_{u\in N(v)}f(uv).$$
The labeling $f$ is called $\Gamma$\textit{-irregular} if for all $v,u\in V(G)$ $w(v)\neq w(u)$.
Furthermore the \textit{group irregularity strength} of $G$, which we denote as $s_g(G)$, is the smallest integer $s$ such that for every Abelian group $\Gamma$ of order $s$ there exists a $\Gamma$-irregular labeling $f$ of $G$.
If we consider the Abelian group to be only a cyclic group $\mathrm{Z}_k$, then we have a \textit{modular edge-graceful labeling} which was introduced by Jones in \cite{phdJones,ModEdGr,Now0ModEdGr}.
He described \textit{modular edge-gracefulness of graphs} as the smallest integer $k(G)= k\geq n$ for which there exists an edge labeling $f:E(G) \rightarrow \mathbb{Z}_k$ such that the induced vertex labeling $g:V(G) \rightarrow \mathrm{Z}_k$ defined as 
$$g(u)=\sum\limits_{v\in N(u)} f(uv) \text{ mod }k$$
is one-to-one.\\

We can easily deduct the correlation between  $k(G)$ and $s_g(G)$. Namely, for every graph $G$ with no component of the order less than 3, there is $k(G) \leq s_g(G)$, \cite{GrIrLaDisc}.

In a paper \cite{GrIrStrCon} Anholcer, Cichacz and Milani$\check{c}$ have shown the following theorem that describes the value of $s_g(G)$ for all connected graphs $G$ of order $n\geq 3$. 
\begin{theorem}[\cite{GrIrStrCon}]\label{Anh}
    Let $G$ be an arbitrary connected graph of order $n\geq 3$. Then
    $$s_g(G) = \left\{\begin{matrix}
        n+2 & \text{if } G\cong K_{1,3^{2q+1}-2} \text{ for some integer } q\geq 1\\ 
        n+1 & \text{if } n\equiv 2 (\text{mod} 4) \wedge G\ncong K_{1,3^{2q+1}-2} \text{ for any integer } q\geq 1 \\ 
        n & \text{otherwise} 
    \end{matrix}\right.$$
\end{theorem}
It is easy to see that to distinguish all $n$ vertices in an arbitrary graph of that order we need at least $n$ different element of $\Gamma$.
Although the following lemma shows that an Abelian group of order $n$ is not always enough to have a $\Gamma$-irregular labeling of $G$.
\begin{lemma}[\cite{GrIrLaDisc}]\label{2mod4}
    Let $G$ be a graph of order $n$.
    If $n\equiv 2 \pmod4$, then there is no $\Gamma$-irregular labeling $G$ for any Abelian group $\Gamma$ of order $n$.
\end{lemma}
So far the best upper bound for any graph was given by Anholcer, Cichacz, and Przybyło in \cite{LinearBoundNo0}
\begin{corollary}[\cite{LinearBoundNo0}]\label{linear}
    Let $G$ be an arbitrary graph of order $n$ having no component of the order less than 3. 
    Then  $s_g(G) \leq 2n.$
\end{corollary}
The following theorem \cite{GrIrLaDisc} by Anholcer and Cichacz determines the value of $k(G)$ for an arbitrary graph $G$ without any components of the order less than 3 or isomorphic to a star of even order.

\begin{theorem}[\cite{GrIrLaDisc}]\label{cykliczne}
    Let $G$ be a graph of order $n$ having neither a component of the order less than 3 nor a $K_{1,2u+1}$ component for any integer $u\geq 1$.
    Then
     $$\begin{matrix}
       k(G) = n & \text{if} & n\equiv 1 \pmod2\\ 
       k(G) = n+1 & \text{if} & n\equiv 2 \pmod4\\
       k(G) \leq n+1 & \text{if} & n\equiv 0 \pmod4.
    \end{matrix}$$
\end{theorem}

In the literature is also considered irregular labeling of digraphs. Cichacz and Tuza showed that if  $n$ is large enough
with respect to an arbitrarily fixed $\varepsilon > 0$ then $\overrightarrow{G}$
has a $\Gamma$-irregular labeling for any $\Gamma$ such that $|\Gamma|>(1+\varepsilon)n$ \cite{CicTuz}. This bound was recently improved for digraphs with no weakly connected components of cardinality less than\/ $4$. Namely, it was shown that they have a $\Gamma$-irregular labeling for every $\Gamma$ such that\/
   $|\Gamma|\geq n+6$ \cite{Cic}. Therefore on one hand we could expect that there exists a constant $K$ such that any graph of order $n$ without components of the order less than $3$ has a $\Gamma$-irregular labeling for any group $\Gamma$ of order at least $n+K$. On the other hand, we see in Theorem~\ref{Anh} that some even stars can cause problems.

In this paper we will generalize Theorem~\ref{cykliczne} for all Abelian groups, hence we show that such a constant exists (and is equal to 1) if a graph $G$ does not have even stars as components. Moreover, we will improve the bound given in Corollary~\ref{linear}. We will show that if $\Gamma$ has order at least  $n+K(q_0)$, where $q_0$ is the number of components being even stars in $G$, then $G$ has a $\Gamma$-irregular labeling. However, this bound is not better than those given in Corollary~\ref{linear} for $G$ being a union of $K_{1,3}$, it is much better for graphs such that the number of even stars is small compared to the order of $G$. In particular, we proved, that if  $G$ is a graph of order $n$ having neither a component of the order less than 3 nor a component being a star $K_{1,1+2u}$ for $u\in\{1,2\}$, then $s(G)\leq n+3$.

\section{Proof of the main result}

In the proof of the main result of this paper, we will use the construction of partition of $\Gamma\setminus\{0\}$ into triplets and pairs, which was introduced in \cite{PartAb}.
We call triplets $T_i = \{a_i,b_i,c_i\}$ and $T_j= \{a_j,b_j,c_j\}$ \textit{complementary} if $-a_i = a_j, -b_i = b_j$ and $-c_i=c_j$.
We will denote the $T_j$ triplet as $-T_i$.
Let $\Gamma$ be an Abelian group of odd order $n=6m+s$ $(s = 1,3 \text{ or } 5)$.
A partition of $\Gamma\backslash \{0\}$ into $m$ pairs of complementary triplets $T_1,-T_1,\ldots, T_m, -T_m$ and $\frac{s-1}{2}$ sets $\{d_j, -d_j\}$, where $\sum\limits_{a\in T_i} a =0$, is called a \textit{Skolem partition} of $\backslash \{0\}$. 
\begin{theorem}[\cite{PartAb}]\label{Skolem}
    Every Abelian group $\Gamma$ of odd order admits a Skolem partition.
\end{theorem}

For any two given vertices $v_1$ and $v_2$ from the same connected component of a graph $G$ exist walks from $v_1$ to $v_2$.
If a walk has an even number of vertices then we are going to call it an \textit{even walk}.
Otherwise, it is called an \textit{odd walk}.
In this proof, we will consider the shortest odd or even walk.\\
While labeling edges of a graph $G$ we will start with $0$ on all of them.
Next, we will choose $v_1$ and $v_2$ and modify all of the edges of a chosen walk from  $v_1$ to $v_2$ by adding some element of an Abelian group $\Gamma$.
Roughly speaking we will be adding some element $a\in \Gamma$ to all the labels of the edges on an odd position on the walk, starting from $v_1$.
However, to the rest of the labels of edges on the walk we will add $-a$. 
This modification of labels we will denote with $\phi_e(v_1, v_2)=a$ if we are considering the shortest even walk and $\phi_o(v_1, v_2)=a$ otherwise.
It is easy to see that $\phi_e(v_1, v_2)=a$ adds $a$ to a weighted degree of both $v_1$ and $v_2$, while $\phi_o(v_1, v_2)=a$ increases the weighted degree of $v_1$ by $a$ and decreases the weighted degree of $v_2$ by $a$.
We should also note that in both cases the operation does not modify the weighted degree of any other vertices of the walk.

We will need the following lemma:
\begin{lemma}\label{lemat}
    Let $G$ be a graph of order $n$ having neither a component of the order less than 3 nor a $K_{1,2u+1}$, $u\geq 1$. Let $\Gamma$ be an Abelian group of order $t\geq n$.
    Let $q_1$ be the number of bipartite components of $G$ with both color classes odd and  $q_2$ be the number of components of odd order. Let $S$ be a subsets of elements of $\Gamma$ such that $S$ contains $2m=2q_1+2\left\lfloor\frac{q_2}{2}\right\rfloor$ complementary triples $T_1,T_2,\ldots, T_{2m-1}, T_{2m}$ such that $T_{2i}=-T_{2i-1}$ for $i=1,2,\ldots,m$ and $k=\left\lfloor\frac{n-2m}{2}\right\rfloor$ zero-sum pairs $D_1,D_2,\ldots,D_k$, and the identity element $\{0\}$.

    Then  there exists a $\Gamma$-irregular labeling $f\colon E(G) \to \Gamma$  such that $w(v) \in S$ for any $v\in V(G)$.
\end{lemma}
\begin{Proof}
    This proof is based on the proof of Lemma 2.5 \cite{GrIrLaDisc}.\\
    Let us denote the elements of a triplet $T_i$ by $x_i, y_i$ and $z_i$ for $i=1,2,\ldots,2m$.
  If $q_2$ is odd, then let $T_{2m+1}=\{0,d_k,-d_k\}$. 
    If a component is a bipartite graph with both color classes even then we divide the vertices of every class into pairs $(v_i^1, v_i^2)$.
    Next we put $$\phi_o(v_i^1, v_i^2) = d_i$$ for every such pair. Therefore $w(v_i^1)=d_i$ and $w(v_i^2)=-d_i$. 
    If a component is not bipartite but is even, then we label it analogically.
    When a component is a bipartite graph with both color classes odd, then we can be sure that there are at least 3 vertices in both classes.
    Let us choose three of them from each class and denote them with $a_i, b_i $ and $c_{i+1}$ in one of the classes and in the second with $a_{i+1}, b_{i+1} $ and $c_i$.
    Than we can put $$\phi_e(a_i,c_i) = y_i,$$
    $$\phi_e(b_i,c_i) = z_i,$$
    $$\phi_e(a_{i+1},c_{i+1}) = y_{i+1},$$
    $$\phi_e(b_{i+1},c_{i+1}) = z_{i+1}.$$
    Thus $w(a_{2j-1})=y_{2j-1}$, $w(b_{2j-1})=z_{2j-1}$, $w(c_{2j-1})=-(y_{2j-1}+z_{2j-1})=-x_{2j-1}$, $w(a_{2j})=y_{2j}=-y_{2j-1}$, $w(b_{2j})=z_{2j}=-z_{2j-1}$ and  $w(c_{2j+1})=x_{2j-1}$ for $j=1,2,\ldots,q_1$. 
    The rest of the vertices we label in the same way as a bipartite component of both classes even.
    If a component is of odd order we choose three vertices $a_i, b_i $ and $c_i$, when it is a bipartite component then $c_i$ we choose from a color class of odd order, and $a_i, b_i $ from the even class. 
    Next we put $$\phi_e(a_i,c_i) = y_i,$$
    $$\phi_e(b_i,c_i) = z_i.$$
    The rest of the vertices we label analogically to earlier cases.\\
    It is easy to see that the weights of vertices are pairwise distinct and $w(v)\in S$ for any $v\in V(G)$.
\end{Proof}

Note that now we easily obtain the following results:
\begin{corollary}
  Let $G$ be a graph of order $n$ having neither a component of the order less than 3 nor a $K_{1,2u+1}$, $u\geq 1$. Then for every odd integer $t\geq n$ there exists a $\Gamma$-irregular labeling, where $\Gamma$ is an Abelian group of order $t$.
\end{corollary}
\begin{Proof}
Let $q_1$ as the number of bipartite components of $G$ with both color classes odd and  $q_2$ be the number of components of odd order. Let  $2m=2q_1+2\left\lfloor\frac{q_2}{2}\right\rfloor$ and $k=\left\lfloor\frac{n-2m}{2}\right\rfloor$.

There exists a Skolem partition of $\Gamma$ by Theorem~\ref{Skolem}. It is easy to see that if a Skolem partition of $\Gamma$ exists, then there is a partition of $\Gamma\setminus\{0\}$ into complementary triples $T_1,T_2,\ldots, T_{2m-1}, T_{2m}$ such that $T_{2i}=-T_{2i-1}$ for $i=1,2,\ldots,m$ and $k'=\left\lfloor\frac{t-2m}{2}\right\rfloor\geq k$ zero-sum pairs $D_1,D_2,\ldots,D_{k'}$.
Now we apply Lemma~\ref{lemat} and we are done.
\end{Proof}
The above corollary and Lemma~\ref{2mod4} imply the following bound analogous to those from Theorem~\ref{cykliczne} for cyclic groups.

\begin{theorem}\label{bez_gwiazd}
    Let $G$ be a graph of order $n$ having neither a component of the order less than 3 nor a $K_{1,2u+1}$ component for any integer $u\geq 1$.
    Then
     $$\begin{matrix}
       s_g(G) = n & \text{if} & n\equiv 1 \pmod2\\ 
       s_g(G) = n+1 & \text{if} & n\equiv 2 \pmod4\\
       s_g(G) \leq n+1 & \text{if} & n\equiv 0 \pmod4.
    \end{matrix}$$
\end{theorem}

For graphs with star components, we will need now the following.

\begin{lemma}\label{lemat_glowny}
    Let $G$ be a graph of order $n$ having no component of the order less than 3.
    Let $q_0$ be the number of components being stars of even order, $q_1$ be the number of bipartite components of $G$ with both color classes odd and  $q_2$ be the number of components of odd order. Let $q_4=\left\lfloor\frac{ n-6\left\lceil\frac{q_0}{2}\right\rceil-6q_1-6\left\lfloor\frac{q_2}{2}\right\rfloor-1}{2}\right\rfloor$ if $q_0$ and $n$ are both odd and $q_4=\left\lfloor\frac{ n-3q_0-6q_1-6\left\lfloor\frac{q_2}{2}\right\rfloor}{2}\right\rfloor$ otherwise. Let  $\varepsilon(q_0)=3$ for $q_0$ odd and $n$ even,  and $\varepsilon(q_0)=0$ otherwise. Let  $\Gamma$ be an Abelian group of odd order $t\geq n+K(q_0,q_4)$, where $$K(q_0,q_4)=\max\left\{ \varepsilon(q_0),5q_0-2q_4+\varepsilon(q_0)-1\right\}.$$
    
    Then  there exists a $\Gamma$-irregular labeling $f\colon E(G) \to \Gamma$.
\end{lemma}
\begin{Proof}
Let $m=3q_1+3\left\lfloor\frac{q_2}{2}\right\rfloor$,  $m'=\left\lceil\frac{q_0}{2}\right\rceil$. Let $k'=\left\lfloor\frac{ t-2m-6m'-1}{2}\right\rfloor$ if $q_0$ and $n$ are both odd, $k'=\left\lfloor\frac{t-2m-6m'}{2}\right\rfloor\geq q_4$ otherwise.
Observe that  $k'\geq 5m'-2$ for $q_0$ even and $k'\geq 5m'-4$ otherwise.

A Skolem partition of $\Gamma$ by Theorem~\ref{Skolem} exists. It is easy to see that if a Skolem partition of $\Gamma$ exists, then there is a partition of $\Gamma\setminus\{0\}$ into complementary triples $T_1,T_2,\ldots, T_{2m-1}, T_{2m}, T_{2m+1}, \ldots,$ $T_{2m+2m'-1}, T_{2m+2m'},$ (i.e. $T_{2i}=-T_{2i-1}$ for $i=1,2,\ldots,m+m'$) and zero-sum pairs $D_1,D_2,\ldots,D_{k'}$.
Let us denote the elements of a triplet $T_i$ by $x_i, y_i$ and $z_i$ for $i=1,2,\ldots,2m+2m'$, and the elements of a pair $D_j$ by $d_j,-d_j$ for $j=1,2,\ldots,k'$. 

Let $K_{1,1+2u_1}$, $K_{1,1+2u_2}$,$\ldots$,$K_{1,1+2u_{q_0}}$ be stars of even order such that $c_i$ is the central vertex and $v_i^1,v_i^2,\ldots,v_i^{1+2u_i}$ are the leaves in $K_{1,1+2u_i}$. Let 
$$\phi_e(c_i,v_i^1) = x_{2m+i},$$
$$\phi_e(c_i,v_i^2) = y_{2m+i},$$
The vertices $v_i^1$ and $v_i^2$ get the permanent weights $w(v_i^1)=x_{2m+i}$, $w(v_i^2)=y_{2m+i}$ resp., whereas we will change the weight of vertices $c_i$ for $i=1,2,\ldots,q_0$.
  
  Suppose first that $q_0$  or $n$ is even.
  Let
$$S_1=\bigcup_{i=1}^{k'}\{d_i,-d_i\}\cup\bigcup_{i=0}^{m'-1}\{z_{2m+i},-z_{2m+i}\},$$
$$W_1=\bigcup_{i=0}^{m'-1}\{x_{2m+2i+1},-x_{2m+2i+1},y_{2m+2i+1},-y_{2m+2i+1}\}.$$
Note that $w(v_i^1),w(v_i^2)\in W_1$. During the process of labeling in each step the set $S_i$ will be the set with available labels and the set $W_i$ is the set of occupied labels (i.e. the permanent weights). In the $(2i+1)$th step for for $i=0,1,\ldots,m'-1$  we will pick an element $g \in S_{2i+1}$ in such a way that $g-z_{2m+i+1}\not\in W_{2i+1}\cup\{0,-g\}$. Note that $|W_{2i+1}|=4m'+4i$ and $|S_{2i+1}|=2k'+2m'-4i$.
Since $\Gamma$ has an odd order there exists exactly one element $x\in \Gamma$ such that $x-z_{2m+i+1}=-x$. Moreover because $k'> 5m'-3$, there is $|W_{2i+1}|+2<|S_{2i+1}|$ and we can always choose such element $g\in S_{2i+1}$ for $i=0,1,\ldots,m'-1$. Then we put 
$\phi_e(c_i,v_i^2) =g$ and $\phi_e(c_{i+1},v_{i+1}^2) =-g$ and we change the set of available labels
$$S_{2(i+1)+1}=S_{2i+1}\setminus\{g,-g,g-z_{2m+i+1},-g+z_{2m+i+1}\}$$ and  the set of permanent weights
$$W_{2(i+1)+1}=W_{2i+1}\cup\{g,-g,g-z_{2m+i+1},-g+z_{2m+i+1}\}.$$
We got  the permanent weights $w(c_{i+1})=g-z_{2m+i+1}$, $w(v_{i+1}^3)=g$ and $w(c_{i+2})=-g+z_{2m+i+1}$, $w(v_{i+2}^3)=-g$ for $i=0,1,\ldots,m'-1$.  One can easily see that $w(c_{i+2}),w(v_{i+2}^3) \not\in W_{2i+1}\cup\{0,g\}$. 
 
Observe that the set $S_{2m'-1}$ contains at least $t-2m-8m'$ elements which form zero-sum pairs. 
The remaining elements of stars $K_{1,1+2u_i}$ we divide  into pairs $(v_i^{2j}, v_i^{2j+1})$ for $j=2,3,\ldots,u_i$.
    Next we put $$\phi_o(v_i^{2j}, v_i^{2j+1}) = g\in S_{2m'-1} $$ for every such pair. Therefore $w(v_i^1)=g$ and $w(v_i^2)=-g$ from the set $S_{2m'-1}$. Finally, we are left with the set that has at least $n-2m-2q_0-\sum_{i=1}^{q_0}2u_i$ elements that form zero-sum pairs, moreover, we have $2m$ complementary triplets $T_1,T_2,\ldots, T_{2m-1}, T_{2m}$. Thus for the remaining components of the graph $G$ we apply Lemma~\ref{lemat} and we are done.

Suppose now that $q_0$ is odd. 
For a star $K_{1,1+2u_1}$ put $\phi_e(c_1,v_1^1) =- x_{2m+1},$
$\phi_e(c_1,v_1^2) = y_{2m+1},$ and $\phi_e(c_1,v_1^3) = 0$. 
 Let
$$S_3=\bigcup_{i=1}^{k'}\{d_i,-d_i\}\cup\bigcup_{i=1}^{m'-1}\{z_{2m+i},-z_{2m+i}\},$$
$$W_3=\bigcup_{i=0}^{m'-1}\{x_{2m+2i+1},-x_{2m+2i+1},y_{2m+2i+1},-y_{2m+2i+1}\}\cup \{z_{2m+1},-z_{2m+1}\}.$$
If additionally, $n$ is odd, then $q_2$ is odd. Let $C$ be a component of the odd order. As in the proof of Lemma~\ref{lemat}, we can label edges of $C$ in such a way that there exist three vertices with weights $x_{2m+1},y_{2m+1}$ and $-z_{2m+1}$, whereas all the remaining vertices in $C$ we can divide into pairs that each pair is connected by an even walk.
Since now $4m'+4i=|W_{2i+1}|+2<|S_{2i+1}|=2k'+2m'-4i+2$ for $i=0,1,\ldots,m'-1$, thus we proceed now analogously as above and we are done.
\end{Proof}

From the above lemma, we obtain immediately the following:

\begin{theorem}\label{bez_K13}
    Let $G$ be a graph of order $n$ having neither a component of the order less than 3 nor a $K_{1,1+2u}$  component for $u\in\{1,2\}$.
    Then
     $$\begin{matrix}
       s_g(G) = n & \text{if} & n\equiv 1 \pmod2\\ 
       s_g(G) \leq n+3 & \text{if} & n\equiv 0 \pmod2.
    \end{matrix}$$
\end{theorem}
\begin{Proof}
If there is no  $K_{1,1+2u}$ for $u\geq 1$, then we are done by Theorem~\ref{bez_gwiazd}. Therefore we can assume that $q_0>0$ and apply Lemma~\ref{lemat_glowny}. Since $u\geq 2$, there is $n\geq 8q_0+6q_1+3q_2$, what implies that $2q_4\geq 5q_0-1$, thus
$K(q_0,q_4)=\varepsilon(q_0)$. If $n$ is odd, then $\varepsilon(q_0)=0$ and $\varepsilon(q_0)\leq 3$ for $n$ even.
\end{Proof}

\begin{figure}[ht!]
\begin{center}
\begin{tikzpicture}[scale=1,style=thick,x=1cm,y=1cm]
\def\vr{3pt} 


\path (0,2) coordinate (a);
\path (1.5,0) coordinate (b);
\path (1.5,2) coordinate (c);
\path (1.5,4) coordinate (d);

\path (2,2) coordinate (a1);
\path (3.5,0) coordinate (b1);
\path (3.5,2) coordinate (c1);
\path (3.5,4) coordinate (d1);

\path (5,2) coordinate (a2);
\path (6.5,0) coordinate (b2);
\path (6.5,2) coordinate (c2);
\path (6.5,4) coordinate (d2);

\path (7.5,2) coordinate (a3);
\path (9,0) coordinate (b3);
\path (9,2) coordinate (c3);
\path (9,4) coordinate (d3);

\path (10.5,2) coordinate (a4);
\path (12,0) coordinate (b4);
\path (12,2) coordinate (c4);
\path (12,4) coordinate (d4);

\path (13.25,2) coordinate (a5);
\path (14.75,0) coordinate (b5);
\path (14.75,2) coordinate (c5);
\path (14.75,4) coordinate (d5);

\draw (b) -- (a) -- (c);
\draw (a) -- (d);

\draw (b1) -- (a1) -- (c1);
\draw (a1) -- (d1);

\draw (b2) -- (a2) -- (c2);
\draw (a2) -- (d2);

\draw (b3) -- (a3) -- (c3);
\draw (a3) -- (d3);

\draw (b4) -- (a4) -- (c4);
\draw (a4) -- (d4);

\draw (b5) -- (a5) -- (c5);
\draw (a5) -- (d5);

\draw (a) [fill=white] circle (\vr);
\draw (b) [fill=white] circle (\vr);
\draw (c) [fill=white] circle (\vr);
\draw (d) [fill=white] circle (\vr);
\draw (a1) [fill=white] circle (\vr);
\draw (b1) [fill=white] circle (\vr);
\draw (c1) [fill=white] circle (\vr);
\draw (d1) [fill=white] circle (\vr);
\draw (a2) [fill=white] circle (\vr);
\draw (b2) [fill=white] circle (\vr);
\draw (c2) [fill=white] circle (\vr);
\draw (d2) [fill=white] circle (\vr);
\draw (a3) [fill=white] circle (\vr);
\draw (b3) [fill=white] circle (\vr);
\draw (c3) [fill=white] circle (\vr);
\draw (d3) [fill=white] circle (\vr);
\draw (a4) [fill=white] circle (\vr);
\draw (b4) [fill=white] circle (\vr);
\draw (c4) [fill=white] circle (\vr);
\draw (d4) [fill=white] circle (\vr);
\draw (a5) [fill=white] circle (\vr);
\draw (b5) [fill=white] circle (\vr);
\draw (c5) [fill=white] circle (\vr);
\draw (d5) [fill=white] circle (\vr);

\draw(0,2.25) node {$4$};
\draw[anchor = east] (b) node {$0$};
\draw[anchor = south] (c) node {$1$};
\draw[anchor = east] (d) node {$3$};

\draw(2.0,2.25) node {$6$};
\draw[anchor = east] (b1) node {$2$};
\draw[anchor = south] (c1) node {$5$};
\draw[anchor = east] (d1) node {$7$};

\draw(4.65,2.25) node {$(2,2)$};
\draw[anchor = east] (b2) node {$(1,1)$};
\draw[anchor = north] (c2) node {$(1,2)$};
\draw[anchor = east] (d2) node {$(0,3)$};

\draw(7.15,2.25) node {$(1,3)$};
\draw[anchor = east] (b3) node {$(1,0)$};
\draw[anchor = north] (c3) node {$(0,1)$};
\draw[anchor = east] (d3) node {$(0,2)$};

\draw(10.15,2.25) node {$(0,0,0)$};
\draw[anchor = east] (b4) node {$(1,1,0)$};
\draw[anchor = north] (c4) node {$(0,1,1)$};
\draw[anchor = east] (d4) node {$(1,0,1)$};

\draw(12.8,2.25) node {$(1,1,1)$};
\draw[anchor = east] (b5) node {$(1,0,0)$};
\draw[anchor = north] (c5) node {$(0,0,1)$};
\draw[anchor = east] (d5) node {$(0,1,0)$};

\draw(1,1) node {$0$};
\draw(0.75,2.25) node {$1$};
\draw(0.75,3.35) node {$3$};

\draw(3,1) node {$2$};
\draw(2.75,2.25) node {$5$};
\draw(2.75,3.35) node {$7$};

\draw(5.25,1) node {$(1,1)$};
\draw(5.7,2.25) node {$(1,2)$};
\draw(5.5,3.35) node {$(0,3)$};

\draw(7.75,1) node {$(1,0)$};
\draw(8.25,2.25) node {$(0,1)$};
\draw(8.05,3.35) node {$(0,2)$};

\draw(10.5,1) node {$(1,1,0)$};
\draw(11.45,2.25) node {$(0,1,1)$};
\draw(10.75,3.35) node {$(1,0,1)$};

\draw(13.35,1) node {$(1,0,0)$};
\draw(14.15,2.25) node {$(0,0,1)$};
\draw(13.65,3.35) node {$(0,1,0)$};
\end{tikzpicture}
\end{center}
\caption{Irregular labelings of $2K_{1,3}$ in groups $\mathbb{Z}_8$, $\mathbb{Z}_2\times \mathbb{Z}_4$ and $\mathbb{Z}_2\times \mathbb{Z}_2\times \mathbb{Z}_2$, respectively.}\label{rys}
\label{less}
\end{figure}
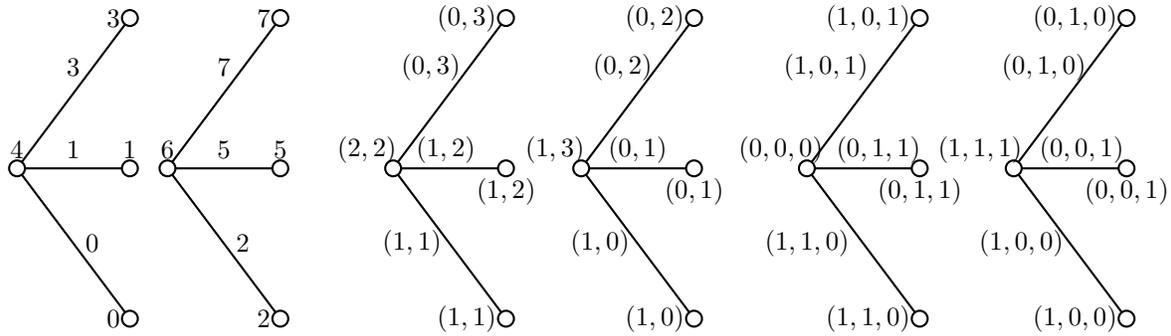

Note that for some families of graphs containing small stars as components we can obtain better a bound than those in Theorem~\ref{linear}. Let $tG$ denote the union of $t$ disjoint copies of $G$, then if $G=tK_{1,3}\cup lC_4$ and $l\geq t$ we have    $s(G)\leq 4(t+l)+3$ by Lemma~\ref{lemat_glowny}. Moreover in Figure~\ref{rys} we show that $s_g(2K_{1,3})=8$. Therefore at the end of this section, we state the following conjecture.
\begin{conj}There exists a constant $K$ such that any graph $G$ of order $n$ with no  connected components
 of order less than $3$ has a $\Gamma$-irregular labeling for every\/ $\Gamma$ such that $|\Gamma|\geq n+K$.
\end{conj}

\bibliographystyle{siam}

\end{document}